\documentclass[11pt,leqno,a4paper,dvips,twoside]{article}

\usepackage{graphicx}
\usepackage{amsfonts}
\usepackage{amssymb}
\usepackage{amsmath}

\newtheorem{theo}{Theorem}%[section]
\newtheorem{prop}{Proposition}%[section]
\newtheorem{rema}{Remark}%[section]
\begin{document}
\thispagestyle{plain} {\footnotesize {{\bf General Mathematics Vol.
22, No. 1 (2014),  19--26}}} \vspace*{2cm}

\begin{center}
%%% TITLE  %%%
{\Large {\bf On Bullen's and related inequalities}} \vspace{0.5cm}
\footnote{\it Received 15 June, 2014

\hspace{0.1cm} Accepted for publication (in revised form) 15 July,
2014}\vspace{0.3cm}

%%%% AUTOR %%%%
{\large Ana Maria Acu, Heiner Gonska}
\end{center}
\vspace{0.5cm}
\begin{center}
 {\small\it Dedicated to the late Academician Professor Dr. Dimitrie
 D. Stancu}
 \end{center}

\begin{abstract}
%%%%% ABSTRACT %%%%%
The estimate in Bullen's inequality will be extended for continuous
functions using the second order modulus of smoothness. A different
form of this inequa-\linebreak lity will be given in terms of the least concave
majorant. Also, the composite case of Bullen's inequality is
considered.

\end{abstract}

\begin{center}

\end{center}
\vspace*{0.3cm}

% Next you must introduce the contents of your article

\section{Motivation}

Over the years it happened during several editions of RoGer - the
 Romanian-German Seminar on Approximation Theory - that
the  second author learned about inequa-\linebreak lities the validity of which
was known for regular (i.e., differentiable) functions only. Using
tools from Approximation Theory, we showed in \cite{Acu&Gonska_2009}
and \cite{Acu&Gonska&Rasa_2011} that such restriction can sometimes
be dropped and that the estimates can be extended to (at least)
continuous functions on the given compact intervals, for example.

Our more general estimates in \cite{Acu&Gonska_2009} and
\cite{Acu&Gonska&Rasa_2011}
 were given in terms of the least concave majorant of the usual first order modulus of
  continuity. Already this is rather a complicated quantity. There we focussed on Ostrowski- and Gr\"uss-type directions.
The best way seems to be that via a certain K-functional. This road
was recently and thoroughly described in a paper by P\u alt\u anea
\cite{Pal}. But the knowledge about this method is much older. See
papers by Peetre \cite{Pee}, Mitjagin and Semenov \cite{MiSe}  as
well as the diploma thesis of Sperling \cite{Spe}, for example.

In Sections 2 and 3 we will consider classical and composite Bullen
functionals. It will become clear there that it is quite natural to
use the second order modulus of smoothness of a continuous function
and a related K-functional.

At the end of this note we will return to the concave majorant and
its significance in the field of Inequalities. For the composite
case and continuously differentiable functions Section 4 contains a
very precise estimate in terms of the majorant.

\section{Bullen's inequality revisited }

This paper is mostly meant to be a contribution to the ever-lasting
discussion on an inequality given by Bullen in \cite{Bul78} (see
also an earlier paper by Hammer \cite{Ham}) in the following form:

\noindent {\bf Theorem A}. {\it If $f$ is convex and integrable,
then
$$
\left( \int_{-1}^1f \right) - 2 f(0) \le f(-1) + f(1) - \int_{-1}^1
f .
$$
If transformed to an arbitrary compact interval $[a,b] \subset
\mathbb{R}$, $a < b$, the equivalent form of the inequality reads
$$
\frac{2}{b-a} \int_a^b f(x) dx \le \frac{f(a) + f(b)}{2} + f \left(
\frac{a+b}{2}\right).
$$}
This is the form we learned about at "RoGer 2014 - Sibiu". In the
talk of Petric\u a Dicu we also learned that in 2000 Dragomir and
Pearce \cite{DrPe} had given the following inequality for functions
$f$ in $C^2 [a,b]$ with known bounds for $f''$:

\begin{theo}\label{t1}
Let $f : [a,b] \to R$ be a twice differentiable function for which
there exist real constants $m$ and $M$ such that
$$
m \le f'' (x) \le M , \quad \mbox{for all}\quad x \in [a,b].
$$
Then
\begin{equation}\label{e1}
\begin{array}{lcl}
\displaystyle m \frac{(b-a)^2}{24} & \le & \displaystyle \frac{f(a) + f(b)}{2} + f \left( \frac{a+b}{2} \right) - \frac{2}{b-a} \int_a^b f(x) dx \\
& \le & \displaystyle M \frac{(b-a)^2}{24} .
\end{array}
\end{equation}
\end{theo}

If we define the Bullen functional $B$ by
$$
B(f) := \frac{f(a)+f(b)}{2} + f \left( \frac{a+b}{2} \right)  -
\frac{2}{b-a} \int_a^b f(x) dx,
$$
we note that $B$ is defined for all functions in $C[a,b]$ and that -
so far - we have the following

\begin{prop}
The Bullen functional $B : C[a,b] \to \mathbb{R}$ satisfies
\begin{itemize}
\item[(i)] $|B(f)| \le 4 ||f||_\infty$ for all $f \in C[a,b]$, $||\cdot ||_\infty$ indicating the sup norm on $[a,b]$.
\item[(ii)] $|B(g)| \le \displaystyle\frac{(b-a)^2}{24} \cdot ||g''||_\infty$ for all $g\in C^2[a,b]$.
\end{itemize}
\end{prop}
We will next explain how to turn this into a more general statement
using the following result from  \cite{GoKo7}:
\begin{theo}\label{t2}
Let $(F,||\cdot||)_F)$ be a Banach space, and let $H : C[a,b] \to
(F,||\cdot||_F)$ be an operator, where
\begin{itemize}
\item[a)] $||H(f+g)||_F \le \gamma (||Hf||_F + ||Hg||_F)$ for all $f,g \in C[a,b]$;
\item[b)] $||Hf||_F \le \alpha \cdot ||f||_{\infty}$ for all $f \in C[a,b]$;
\item[c)] $||Hg||_F \le \beta_0 \cdot ||g||_{\infty} + \beta_1 \cdot ||g'||_{\infty} + \beta_2 \cdot ||g''||_{\infty}$ for all $g \in C^2 [a,b]$.
\end{itemize}

Then for all $f \in C[a,b]$, $0 < h \le (b-a)/2$, the following
inequality holds:
$$
||Hf||_F \!\le \gamma \left\{ \beta_0 \cdot ||f||_{\infty}\! +\!
\frac{2\beta_1}{h} \omega_1 (f;h) \!+\! \frac34 \left( \alpha\! +\!
\beta_0 \!+\! \frac{2\beta_1}{h}\! + \!\frac{2\beta_2}{h^2} \right)
\omega_2 (f;h)\right\} .
$$
\end{theo}

Taking $H = B$ and $F=\mathbb{R}$  in  Theorem \ref{t2} we have the
following list of constants:
$$
\gamma = 1 , \;\; \alpha = 4 , \;\; \beta_0 = 0 , \;\; \beta_1 = 0,
\;\; \beta_2 = \frac{(b-a)^2}{24} .
$$
This takes us to the following
\begin{prop}
For the Bullen functional $B$ and all $f \in C[a,b]$ we have
$$
|B(f)| \le \left( 3 + \left( \frac{b-a}{4h}\right)^2\right) \cdot
\omega_2 (f,h) , \;\; 0 < h \le \frac{b-a}{2} .
$$
The special choice $h =\displaystyle \frac{b-a}{k}$, $k\geq 2$
yields
$$
|B(f)| \le \left(3+\displaystyle\frac{k^2}{16}\right) \cdot \omega_2
\left( f , \frac{b-a}{k}\right) .
$$
\end{prop}

\begin{rema}
For $f \in C^2[a,b]$ the latter inequality implies
$$
|B(f)| \le \left(\displaystyle \frac{1}{16}+\frac{3}{k^2}\right)
(b-a)^2||f''||_{\infty},  k\geq 2.
$$
As far as the constant $ \left(\displaystyle
\frac{1}{16}+\frac{3}{k^2}\right)$ is concerned, this is much worse
than what was invested for $C^2$ functions.
\end{rema}

An alternative estimate is given in the next proposition. Note that it also follows from Theorem 6 in Gavrea's paper \cite{Gavrea}.

\begin{prop}
If the second $K$-functional on $C[a,b]$ is defined  by
$$
%\begin{array}{l}
K(f;t^2;C[a,b],C^2[a,b])% \\[4mm]
\displaystyle := \inf \{||f-g||_\infty + t^2 ||g''||_{\infty} : g
\in C^2[a,b]\} , \;\; t \ge 0,
%\end{array}
$$
then
$$
|B(f)| \le
4K\left(f;\displaystyle\frac{(b-a)^2}{96};C[a,b],C^2[a,b]\right).
$$
\end{prop}

\bigskip

\noindent{\bf Proof.} In \cite{Acu&Gonska_2009} the following result
was obtained:
\begin{align*}
&\left|\displaystyle\frac{1}{2}\left[(x-a)f(a)+(b-a)f(x)+(b-x)f(b)\right]-\displaystyle\int_{a}^{b}f(t)dt\right|\\
&\leq
2(b-a)K\left(f;\displaystyle\frac{(x-a)^3+(b-x)^3}{24(b-a)};C[a,b],C^2[a,b]\right).
\end{align*}
The proposition is proved if we substitute
$x=\displaystyle\frac{a+b}{2}$ in the above  inequality.

\begin{rema}
\begin{itemize}
\item[(i)] For $h \in C^2 [a,b]$, we have
$$
\begin{array}{lcl}
|B(h)| & \le & \displaystyle 4 \cdot K\left( h;\frac{(b-a)^2}{96};C[a,b],C^2[a,b]\right) \\[4mm]
& = & \displaystyle 4 \cdot \inf \left\{ ||h-g||_\infty + \frac{(b-a)^2}{96} ||g''||_{\infty} , g \in C^2[a,b]\right\} \\[4mm]
& \le & \displaystyle 4 \cdot \frac{(b-a)^2}{96} ||h''||_{\infty}\,\textrm{(taking g=h)} \\[4mm]
& = & \displaystyle�\frac{(b-a)^2}{24} ||h''||_{\infty}, \;\;
\mbox{i.e.},
\end{array}
$$
 the original inequality for $C^2$ functions.
\item[(ii)] It is known that, for $f \in C[a,b]$,
$$
K(f;t^2;C[a,b],C^2[a,b]) \le c \cdot \omega_2 (f,t) , \;\; 0 \le t
\le \frac{b-a}{2}
$$
 with a constant $c \not= c(f,t)$. According to our knowledge the best possible value of $c$ is unknown.

Zhuk showed in \cite{Zhu1989} that, for $t \le\displaystyle
\frac{b-a}{2}$, one has
$$
K(f;t^2; C[a,b], C^2[a,b]) \le \frac94 \cdot \omega_2 (f;t).
$$
Using the latter we arrive, for $h\in C^2[a,b]$, at
$$ \left| B(h) \right|\leq \displaystyle\frac{3(b-a)^2}{32}\|h^{\prime\prime}\|_{\infty}. $$
\end{itemize}
\end{rema}
\section{ The composite case}
Here we consider the composite case of the Bullen functional, i.e.,
the functional which arises when comparing the composite trapezoidal
and midpoint rules. To this end the interval $[a,b]$ is divided in
$n \ge 1$ subintervals as follows
$$ a=x_0<\cdots<x_i<x_{i+1}<\cdots<x_n=b. $$
Let the composite Bullen functional $B_c : C[a,b] \to \mathbb{R}$ be
given by
\begin{align*} B_{c}(f)&=\displaystyle\frac{1}{b-a}\sum_{i=0}^{n-1}(x_{i+1}-x_{i})\left[\displaystyle\frac{f(x_i)+f(x_{i+1})}{2}+f\left(\frac{x_i+x_{i+1}}{2}\right)\right.\\
&\left.-\frac{2}{x_{i+1}-x_{i}}\displaystyle\int_{x_{i}}^{x_{i+1}}f(x)dx\right].
\end{align*}
\begin{prop}\label{p4}
In the composite case there holds
\begin{itemize}
\item[(i)] $|B_c(f)| \le 4 ||f||_\infty$ for all $f \in C[a,b]$,
\item[(ii)] $|B_c(g)| \le \displaystyle\frac{1}{24(b-a)}\sum_{i=0}^{n-1}(x_{i+1}-x_{i})^3||g''||_\infty$ for all $g\in C^2[a,b]$.
\end{itemize}
\end{prop}
Using Theorem \ref{t2} and Proposition \ref{p4} we obtain the
following inequality for the composite Bullen functional involving
the second modulus of continuity:
\begin{prop} For the composite Bullen functional one has
$$|B_c(f)|\leq \left(3+\displaystyle\frac{1}{16(b-a)h^2}\sum_{i=0}^{n-1}(x_{i+1}-x_{i})^3\right)\omega_{2}(f,h),\, 0<h\leq\displaystyle\frac{b-a}{2}.  $$
For
$h=\displaystyle\frac{1}{k}\sqrt{\frac{\sum_{i=0}^{n-1}(x_{i+1}-x_{i})^3}{b-a}},\,
k\geq 2$, this yields
$$ |B_{c}(f)|\leq \left( 3+\displaystyle\frac{k^2}{16}\right)\omega_{2}
\left(f,\displaystyle\frac{1}{k}\sqrt{\frac{\sum_{i=0}^{n-1}(x_{i+1}-x_{i})^3}{b-a}}\right),\,
k\geq 2. $$
\end{prop}
\begin{rema} For $f\in C^2[a,b]$ the latter inequality implies
$$ |B_{c}(f)|\!\leq\! \left(\displaystyle\frac{1}{16}\!+\!\frac{3}{k^2}\right)\displaystyle\frac{\sum_{i=0}^{n-1}(x_{i+1}\!-\!x_{i})^3}{b-a}\|f^{\prime\prime}\|_{\infty}
\!\leq\!
\left(\displaystyle\frac{1}{16}\!+\!\frac{3}{k^2}\right)(b-a)^2\|f^{\prime\prime}\|_{\infty}.
$$
The requirement
$F(x_0,x_1,\dots,x_n)=\displaystyle\sum_{i=0}^{n-1}(x_{i+1}-x_{i})^3\to$
minimum entails $x_{i+1}-x_{i}=\displaystyle\frac{b-a}{n}$,
$i=\overline{0,n-1}.$
\end{rema}

The inequality involving a K-functional is given next.
\begin{prop} For the functional $B_{c}$ given as above, $f\in
C[a,b]$, we have
 \begin{equation}\label{I} \left|B_{c}(f)\right|\leq
4
K\left(f;\displaystyle\frac{1}{96(b-a)}\sum_{i=0}^{n-1}(x_{i+1}-x_i)^3;C[a,b],C^2[a,b]\right).
\end{equation}
\end{prop}
\noindent{\bf Proof.} Let $g\in C^2[a,b]$ arbitrary. Then, for
$f\in
 C[a,b]$,
 \begin{align*}|B_{c}(f)|&\leq |B_{c}(f-g)|+|B_{c}(g)|\\
 &\leq 4\|
 f-g\|_{\infty}+\displaystyle\frac{1}{24(b-a)}\sum_{i=0}^{n-1}(x_{i+1}-x_{i})^3\|g^{\prime\prime}\|_{\infty}\\
 &= 4\left\{\|
 f-g\|_{\infty}+\displaystyle\frac{1}{96(b-a)}\sum_{i=0}^{n-1}(x_{i+1}-x_{i})^3\|g^{\prime\prime}\|_{\infty}\right\}.
  \end{align*}
  Passing to the infimum over $g$ yields inequality (\ref{I}).
  \section{Composite Bullen functional for $f\in C^1[a,b]$}
  Using the least concave majorant of the modulus of continuity in
  this section we consider Bullen's inequality for  $f\in
  C^1[a,b]$.
  \begin{prop}\label{p7}
  If $f\in C^1[a,b]$, then
  $$ |B_{c}(f)|\leq\tilde{\omega}\left(f^{\prime},\displaystyle\frac{1}{24(b-a)}\sum_{i=0}^{n-1}(x_{i+1}-x_{i})^3\right).  $$
  \end{prop}
  \noindent{\bf Proof.} We have
  \begin{align}
  |B_{c}(f)|&\leq\displaystyle\frac{1}{b-a}\sum_{i=0}^{n-1}(x_{i+1}-x_{i})\left|\displaystyle\frac{f(x_i)+f(x_{i+1})}{2}+f\left(\frac{x_i+x_{i+1}}{2}\right)\right.\nonumber\\
  &-\left.\frac{2}{x_{i+1}-x_{i}}\int_{x_{i}}^{x_{i+1}}f(x)dx\right|\nonumber\\
  &=\displaystyle\frac{1}{b-a}\sum_{i=0}^{n-1}(x_{i+1}\!\!-x_{i})\left|\displaystyle\frac{f(x_i)+f(x_{i+1})}{2} -
  \frac{1}{x_{i+1}-x_{i}}\int_{x_{i}}^{x_{i+1}}f(x)dx\right.\nonumber\\
 & +f\left(\frac{x_i+x_{i+1}}{2}\right)
  -\left.\frac{1}{x_{i+1}-x_{i}}\int_{x_{i}}^{x_{i+1}\!\!}f(x)dx\right|\nonumber\\
  &=\displaystyle\frac{1}{b-a}\sum_{i=0}^{n-1}\left|
  \int_{x_i}^{x_{i+1}}\left(\frac{f(x_i)-f(x)}{2}\right.\right.
 + \frac{f(x_{i+1})-f(x)}{2}\nonumber\\
 &+\left.\left.f\left(\frac{x_i+x_{i+1}}{2}\right)-f(x)\right)dx\right|\leq 2\|f^{\prime}\|_{\infty}.\nonumber
  \end{align}
Let $g\in C^2[a,b]$. Using Proposition \ref{p4} and the latter
inequality implies
\begin{align*}
\left|B_{c}(f)\right|&=|B_{c}(f-g+g)|\leq |B_{c}(f-g)|+|B_{c}(g)|\\
&\leq
2\|(f-g)^{\prime}\|_{\infty}+\displaystyle\frac{1}{24(b-a)}\sum_{i=0}^{n-1}(x_{i+1}-x_{i})^3\|g^{\prime\prime}\|_{\infty}\\
%\end{align*}
%\begin{align*}
&=2\left\{\|(f-g)^{\prime}\|_{\infty}+\displaystyle\frac{1}{48(b-a)}\sum_{i=0}^{n-1}(x_{i+1}-x_{i})^3\|g^{\prime\prime}\|_{\infty}\right\}.
\end{align*}
Passing to the infimum over $g\in C^2[a,b]$ we have
$$ |B_{c}(f)|\leq 2K\left(f^{\prime};\frac{1}{48(b-a)}\sum_{i=0}^{n-1}(x_{i+1}-x_{i})^3;C^1[a,b],C^2[a,b]\right), $$
so the result follows as a consequence of the relation (see \cite
{Pal})
$$K\left(f^{\prime};t;C^{1}[a,b],C^2[a,b]\right)=\frac{1}{2}\tilde{\omega}(f^{\prime},2t), 0 \le t \le \frac{b-a}{2}.  $$

A particular consequence of Proposition \ref{p7} is the following
version of Bullen's inequality.
\begin{prop} If $f\in C^1[a,b]$, then
$$ |B(f)|\leq\tilde{\omega}\left(f^{\prime},\displaystyle\frac{(b-a)^2}{24}\right) .$$
\end{prop}

\noindent{\bf Acknowledgment.} The authors most gratefully
acknowledge the efficient help of Birgit Dunkel (University of
Duisburg-Essen) and Elsa while preparing this note.

\bigskip
\noindent
 $\begin{array}{ll}
\textrm{\bf Ana Maria Acu}\\
\textrm{Lucian Blaga University of Sibiu}\\
\textrm{Faculty of Sciences}\\
\textrm{Department of Mathematics and Informatics}\\
\textrm{Str. Dr. I. Ratiu, No.5-7, 550012  Sibiu, Romania}\\
\textrm{e-mail: acuana77@yahoo.com}
\end{array}$

\bigskip
\noindent
 $\begin{array}{ll}
\textrm{\bf Heiner Gonska}\\
\textrm{University of Duisburg-Essen}\\
\textrm{Faculty of Mathematics}\\
%\textrm{Department}\\
\textrm{Forsthausweg 2, 47057 Duisburg, Germany}\\
\textrm{e-mail: heiner.gonska@uni-due.de}
\end{array}$

\end{document}